\documentclass[a4paper,12pt]{amsart}
\usepackage{amssymb}
\usepackage{ifthen}
\nonstopmode \numberwithin{equation}{section}
\newtheorem{thm}{Theorem}
\newtheorem{cor}{Corollary}
\newtheorem{lem}{Lemma}

\newtheorem{conj}{Conjecture}

\theoremstyle{definition}
\newtheorem{defn}{Definition}[section]

\newtheorem{prob}[equation]{Problem}

\newenvironment{rem}{%
\bigskip
\noindent \textsl{{\sl Remark. }}}{\bigskip}
\newenvironment{rems}{%
\bigskip
\noindent \textsl{{\sl Remarks. }}}{\bigskip}

\newcounter {own}
\def\theown {\thesection       .\arabic{own}}

\newenvironment{pf}[1][]{%
 \vskip 3mm
 \noindent
 \ifthenelse{\equal{#1}{}}%
  {{\slshape Proof. }}%
  {{\slshape #1.} }%
 }%
{\qed\bigskip}

\newcounter{alphabet}
\newcounter{tmp}
\newenvironment{Thm}[1][]{\refstepcounter{alphabet}%
\bigskip%
\noindent%
{\bf Theorem \Alph{alphabet}}%
\ifthenelse{\equal{#1}{}}{}{ (#1)}%
{\bf .} \itshape}{\vskip 8pt}

\newcommand{\ID}{{\mathbb D}}

\newcommand{\IC}{{\mathbb C}}
\newcommand{\C}{{\mathbb C}}
\newcommand{\sphere}{{\widehat{\mathbb C}}}

\newcommand{\inv}{^{-1}}
\newcommand{\area}{{\operatorname{Area}\,}}

\newcommand{\dz}{{\partial}}
\newcommand{\dzb}{{\bar\partial}}
\newcommand{\aand}{{\quad\text{and}\quad}}

\def\be{\begin{equation}}
\def\ee{\end{equation}}

\newcommand{\bee}{\begin{enumerate}}
\newcommand{\eee}{\end{enumerate}}

\newcommand{\blem}{\begin{lem}}
\newcommand{\elem}{\end{lem}}
\newcommand{\bthm}{\begin{thm}}
\newcommand{\ethm}{\end{thm}}
\newcommand{\bcor}{\begin{cor}}
\newcommand{\ecor}{\end{cor}}
\newcommand{\beg}{\begin{examp}}
\newcommand{\eeg}{\end{examp}}
\newcommand{\begs}{\begin{examples}}
\newcommand{\eegs}{\end{examples}}
\newcommand{\bdefe}{\begin{defn}}
\newcommand{\edefe}{\end{defn}}
\newcommand{\bprob}{\begin{prob}}
\newcommand{\eprob}{\end{prob}}
\newcommand{\bei}{\begin{itemize}}
\newcommand{\eei}{\end{itemize}}

\newcommand{\bcon}{\begin{conj}}
\newcommand{\econ}{\end{conj}}
\newcommand{\bcons}{\begin{conjs}}
\newcommand{\econs}{\end{conjs}}
\newcommand{\bprop}{\begin{propo}}
\newcommand{\eprop}{\end{propo}}
\newcommand{\br}{\begin{rem}}
\newcommand{\er}{\end{rem}}
\newcommand{\brs}{\begin{rems}}
\newcommand{\ers}{\end{rems}}
\newcommand{\bo}{\begin{obser}}
\newcommand{\eo}{\end{obser}}
\newcommand{\bos}{\begin{obsers}}
\newcommand{\eos}{\end{obsers}}
\newcommand{\bpf}{\begin{pf}}
\newcommand{\epf}{\end{pf}}
\newcommand{\ba}{\begin{array}}
\newcommand{\ea}{\end{array}}
\newcommand{\beq}{\begin{eqnarray}}
\newcommand{\beqq}{\begin{eqnarray*}}
\newcommand{\eeq}{\end{eqnarray}}
\newcommand{\eeqq}{\end{eqnarray*}}

\newcommand{\ra}{\rightarrow}

\newcounter{minutes}
\divide\time by 60
\newcounter{hours}
\multiply\time by 60 \addtocounter{minutes}{-\time}

\begin{document}
\bibliographystyle{amsplain}
\title{Quasiconformal extension of meromorphic functions with nonzero pole}


\author{B. Bhowmik${}^{~\mathbf{*}}$}
\address{Bappaditya Bhowmik, Department of Mathematics,
Indian Institute of Technology Kharagpur, Kharagpur - 721302, India.}
\email{bappaditya@maths.iitkgp.ernet.in}
\author{G. Satpati}
\address{Goutam Satpati, Department of Mathematics,
Indian Institute of Technology Kharagpur, Kharagpur - 721302, India.}
\email{goutam.satpati@gmail.com}
\author{T. Sugawa}
\address{Toshiyuki Sugawa, Graduate School of Information Sciences,
Tohoku University, Aoba-ku, Sendai 980-8579, Japan}
\email{sugawa@math.is.tohoku.ac.jp}

\keywords{ Quasiconformal map, Convolution \\
${}^{\mathbf{*}}$ Corresponding author}

\begin{abstract}
In this note, we consider meromorphic univalent functions $f(z)$
in the unit disc with a simple pole at $z=p\in(0,1)$ which have a
$k$-quasiconformal
extension to the extended complex plane $\sphere,$
where $0\leq k < 1$.
We denote the class of such functions by $\Sigma_k(p)$.
We first prove an area theorem for functions in this class.
Next, we derive a sufficient condition for meromorphic functions
in the unit disc with a simple pole at $z=p\in(0,1)$
to belong to the class $\Sigma_k(p)$.
Finally, we
give a convolution property for functions in the class $\Sigma_k(p)$.
\end{abstract}
\thanks{The first author of this article would like to thank
NBHM, DAE, India (Ref.No.- 2/48(20)/2012/NBHM(R.P.)/R\&D II/14916) for its financial support
and the third author would like to thank JSPS Grant-in-Aid for Scientific Research (B) 22340025 for its partial financial support.}

\maketitle
\pagestyle{myheadings}
\markboth{B. Bhowmik, G. Satpati and T. Sugawa}{Quasiconformal extension}

\bigskip
\section{Introduction}
Let $\IC$ denote the complex plane and $\sphere$ denote the extended complex
plane $\C\cup\{\infty\}.$
We shall use the following notations:
$\ID=\{z : |z|<1\}$,
$\overline{\ID}=\{z : |z|\leq 1\}$,
$\ID^{*}= \{z : |z|>1\}$,
$\overline{\ID^{*}}= \{z : |z|\geq1\}$.
Let $f$ be a meromorphic and univalent function in the unit disk $\ID$
with a simple pole at $z=p\in [0,1)$ of residue $1.$
Since $f(z)-1/(z-p)$ is analytic in $|z|<1,$
one has an expression of the form
\be\label{eq0a}
f(z)=\frac{1}{z-p}+\sum_{n=0}^{\infty}a_n z^n
\ee
in $|z|<1.$
We denote the class of such functions by $\Sigma(p)$.
Let $\Sigma^0(p)$ be the subclass of $\Sigma(p)$
consisting of those functions $f$ for which $a_0=0$ in the
above expansion.
Note that if $f, g\in\Sigma^0(p)$ are related by $g=M\circ f$
for a M\"obius transformation $M,$ then $f=g.$

For a given number $0\le k<1,$
$\Sigma_k(p)$ stands for the class of those functions in $\Sigma(p)$
which admit $k$-quasiconformal extension to the extended plane $\sphere.$
Here, a mapping $F:\sphere\to\sphere$ is called $k$-quasiconformal
if $F$ is a homeomorphism and has locally $L^2$-derivatives on
$\IC\setminus\{F\inv(\infty)\}$
(in the sense of distribution) satisfying $|\dzb F|\le k|\dz F|$ a.e.,
where $\dz F=\partial F/\partial z$ and $\dzb F=\partial F/\partial\bar z.$
Note that such an $F$ is called $K$-quasiconformal more often,
where $K=(1+k)/(1-k)\ge1,$ in the literature.
The quantity $\mu=\dzb F/\dz F$ is called the complex dilatation of $F.$
See the standard textbook \cite{LV:qc} by Lehto and Virtanen
for basic properties of quasiconformal mappings.
Set $\Sigma_k^0(p)=\Sigma^0(p)\cap\Sigma_k(p).$

O.~Lehto \cite{Lehto71} refined the Bieberbach-Gronwall area theorem
to the functions in $\Sigma_k(0)$ in the following form.

\begin{Thm}\label{TheoL}
Let $0\le k<1.$
Suppose that $f(z)=z\inv+a_0+a_1z+a_2z^2+\dots$
is a function in $\Sigma_k(0).$
Then
$$
\sum_{n=1}^\infty n|a_n|^2\le k^2.
$$
Here, equality holds if and only if
$$
f(z)=\frac1z+a_0+a_1z,\quad |z|<1,
$$
with $|a_1|=k.$
Moreover, its $k$-quasiconformal extension is given by setting
$$
f(z)=\frac1z+a_0+\frac{a_1}{\bar z} \qquad\text{for}~|z|\ge1.
$$
\end{Thm}

On the other hand, the area theorem was extended by
P.~N.~Chichra \cite{Chi69} to functions in $\Sigma(p)$ as follows.

\begin{Thm}\label{TheoA}
Let $f\in \Sigma(p)$ have the  expansion in \eqref{eq0a}.
Then
\be\label{eq00a}
\sum_{n=1}^{\infty}n|a_n|^2\leq \frac{1}{(1-p^2)^2}.
\ee
Equality holds for the function
$$
f_p(z)=\frac{1}{z-p}+a_0+\frac{z}{1-p^2}.
$$
\end{Thm}

Our first result establishes an area theorem for the class $\Sigma_k(p).$
Interestingly, the form of extremal functions is different from that of
the function $f_p$ in Theorem B.

\bthm\label{Theorem-1}
Let $0\le k<1$ and $0\le p<1.$
Suppose that $f\in \Sigma_k(p)$ is expressed in the form
of \eqref{eq0a}.
Then
\be\label{eq2a0}
\sum_{n=1}^{\infty}n|a_n|^2\leq \frac{k^2}{(1-p^2)^2}.
\ee
Here, equality holds if and only if $f$ is of the form
\beq\label{eq2a1}
f(z) = \frac1{z-p}+a_0+\frac{a_1z}{1-pz}
\qquad\text{for}~|z|<1,
\eeq
where $a_0$ and $a_1$ are constants with $|a_1|=k.$
Moreover, a $k$-quasiconformal extension of this $f$ is given by setting
\be\label{eq:ext}
f(z) = \frac1{z-p}+a_0+\frac{a_1}{\bar z-p}
\qquad\text{for}~|z|\ge1.
\ee
\ethm

Observe that this is a natural extension of Theorem A.
We remark that the function in \eqref{eq2a1} belongs to $\Sigma(p)$
as long as $|a_1|\le1$
(see the latter part of the proof of Theorem \ref{Theorem-1} below).
This function with $|a_1|=1$ provides another extremal case in \eqref{eq00a}.
As an immediate corollary of the theorem, we obtain the following.

\bcor
Let $0<p<1$ and $0<k<1.$
For $f\in \Sigma_k(p)$ with the expansion \eqref{eq0a}, the following
inequality holds:
$$
|a_1|< \frac{k}{1-p^2}.
$$
\ecor

Note that the inequality $|a_1|\le1/(1-p^2)$ for $f\in\Sigma(p)$ is
sharp in view of Theorem B.
We have no exact value of the best upper bound, say, $M(p,k)$
of $|a_1|$ for $f\in\Sigma_k(p).$
The extremal function in Theorem \ref{Theorem-1} and compactness
of the class $\Sigma_k^0(p)$ yield, at least, the estimates
$k\le M(p,k)<k/(1-p^2)$ for $p, k\in(0,1).$

Secondly, we provide a sufficient condition for functions of the form
\eqref{eq0a} to belong to the class $\Sigma_k(p)$.

\bthm\label{Theorem-2}
Let $0\le k<1$ and $0\le p<1.$
Suppose that $\omega$ is an analytic function in $\ID$
such that $|\omega'(z)|\leq k(1+p)^{-2}$ for $z\in\ID$.
Then the function $f$ given by
$$
f(z)=\frac{1}{z-p}+\omega(z), \quad |z|<1,
$$
is a member of $\Sigma_k(p)$.
A $k$-quasiconformal extension is given by setting
\be\label{eq0}
f(z)=\frac{1}{z-p}+\omega(1/\bar z),\quad \text{for $|z|>1$}.
\ee
\ethm

We note that J.~G.~Krzy\.{z} \cite{Kr76} proved this theorem when $p=0.$
He also gave a convolution theorem in the same paper \cite{Kr76}.
We can also extend it to a modified convolution.
The modified Hadamard product (or the modified convolution)
$f\star g$ of two functions $f, g\in \Sigma(p)$
with expansions
\be\label{eq0b}
f(z)=\frac1{z-p}+\sum_{n=0}^\infty a_nz^n
\aand
g(z)=\frac{1}{z-p}+\sum_{n=0}^{\infty}b_n z^n,\quad z\in \ID
\ee
is defined by
\be\label{eq0c}
(f\star g)(z)=\frac{1}{z-p}+\sum_{n=0}^{\infty}a_n\,b_n z^n,\quad z\in \ID.
\ee

Our third result concerns this Hadamard product.

\bthm\label{Theorem-3}
Let $f\in \Sigma_{k_1}(p)$ and $g\in \Sigma_{k_2}(p)$ for some
$k_1, k_2,p\in[0,1).$
If $\alpha=k_1k_2(1-p)^{-2}<1,$ then the modified Hadamard product
$f\star g$ belongs to $\Sigma_{\alpha}(p).$
\ethm

As we mentioned above, this result reduces to a theorem due to Krzy\.{z}
\cite{Kr76} when $p=0.$

We prove Theorem \ref{Theorem-1} in Section 2 and
Theorems \ref{Theorem-2} and \ref{Theorem-3} in Section 3.
We also give another proof of a part of Theorem \ref{Theorem-1}
as a concluding remark in Section 3.

\section{Proof of Theorem \ref{Theorem-1}}
We start the section with proving the area theorem for functions
in the class $\Sigma_k(p)$.
We follow the idea due to Lehto \cite{Lehto71}.

\bpf[Proof of Theorem \ref{Theorem-1}]
Let $f\in\Sigma_k(p)$ have the expansion in \eqref{eq0a}.
We may suppose that $f$ is already extended to a $k$-quasiconformal
mapping of $\sphere$ onto itself.
If $k=0,$ then the assertion clearly holds.
Hence, we assume that $k>0$ in the rest of the proof.
To start with, we first make a change of variables.
We define $\phi:\sphere\to\sphere$ by $\phi(\zeta)=f(1/\zeta).$
Note that $\phi$ has locally $L^2$-derivatives on
$\IC\setminus\{\phi\inv(\infty)\}=\C\setminus\{1/p\}.$
Since the function $\psi(\zeta)=\phi(\zeta)-\zeta/(1-p\zeta)$ has the expression
\be\label{eq:psi}
\psi(\zeta)=\phi(\zeta)-\frac{\zeta}{1-p\zeta}
=\sum_{n=0}^\infty\frac{a_n}{\zeta^n}, \quad |\zeta|>1,
\ee
and, in particular, is bounded and analytic near the point $\zeta=1/p,$
the function $\psi$ has locally $L^2$-derivatives on $\C.$
Therefore for every $r>0$, we can apply the Cauchy-Pompeiu formula
(see \cite[III \S 7]{LV:qc} for details) to the function
$\psi$ in the disk $|\zeta|<r$ to obtain
$$
\psi(\zeta)
= \frac{1}{2\pi i}\int\limits_{|w|=r}\frac{\psi(w)}{w -\zeta}\,dw -
\frac{1}{\pi}{\displaystyle \iint\limits_{|w|<r}}\frac{\dzb\psi(w)}{w-\zeta}\,
dudv,
$$
where $w=u+iv.$
We note that $\psi(\zeta)\to a_0$ as $\zeta\ra \infty$ and
$\dzb\psi(\zeta)=0$ for $|\zeta|>1$.
Letting $r\to+\infty,$ we thus get
\be\label{eq1a}
\psi(\zeta)
=a_0-\frac{1}{\pi}\iint\limits_{|w|<1}\frac{\dzb\psi(w)}{w-\zeta}\,dudv,
\quad \zeta\in \IC.
\ee
We differentiate the above expression with respect to $\zeta$ and obtain
\be\label{eq2a}
\dz\psi(\zeta)
=-\frac{1}{\pi}\iint\limits_{|w|<1}\frac{\dzb\psi(w)}{(w-\zeta)^2}\,dudv
= H[\dzb\psi](\zeta), \quad \zeta\in \IC,
\ee
where $H$ is the two dimensional Hilbert transformation.
(Strictly speaking, the above integral should be understood as Cauchy's
principal value for $|\zeta|\le1.$ See \cite[III \S 7]{LV:qc} for details.)
Since $H$ is a linear isometry of $L^2(\IC),$ in conjunction with \eqref{eq2a},
we have
\be\label{eq2b}
\iint\limits_{\ID}|\dzb\psi(\zeta)|^2\,d\xi d\eta
=\iint\limits_{\IC}\big|H[\dzb\psi](\zeta)\big|^2\,d\xi d\eta
=\iint\limits_{\IC}|\dz\psi(\zeta)|^2\,d\xi d\eta,
\ee
where $\zeta = \xi + i\eta$.
Next, we recall that Chichra indeed showed the following relation in
the proof of Theorem B:
$$
\area(\IC\setminus f(\ID))
=\pi\left[ \frac{1}{(1-p^2)^2}-\sum_{n=1}^{\infty}n|a_n|^2\right].
$$
We remark that $f(\partial\ID)$ is of area zero because $f$ is quasiconformal.
Noting $\phi(\ID)=f(\ID^*)=\IC\setminus f(\overline{\ID}),$
we thus have the relation
\be\label{eq:Chichra}
\area \phi(\ID)
=\pi\left[ \frac{1}{(1-p^2)^2}-\sum_{n=1}^{\infty}n|a_n|^2\right].
\ee
Since $|\dzb \phi|\le k|\dz \phi|$ a.e.,
the Jacobian $J_\phi$ of $\phi$ satisfies the inequality
$$
J_\phi=|\dz \phi|^2-|\dzb \phi|^2\ge(1-k^2)|\dz\phi|^2
\ge(k^{-2}-1)|\dzb\phi|^2
=(k^{-2}-1)|\dzb\psi|^2.
$$
Hence, we obtain
\be\label{eq:area}
\area \phi(\ID)=
\iint\limits_{\ID}J_\phi(\zeta)d\xi d\eta
\ge(k^{-2}-1)\iint\limits_{\ID}|\dzb\psi(\zeta)|^2d\xi d\eta.
\ee
Next, we see from \eqref{eq2b} that
\be\label{eq2d}
\iint\limits_{\ID}|\dzb\psi(\zeta)|^2\,d\xi d\eta=
\iint\limits_{\IC}|\dz\psi(\zeta)|^2\,d\xi d\eta
\geq \iint\limits_{|\zeta|>1}|\dz\psi(\zeta)|^2\,d\xi d\eta.
\ee
It is easy to evaluate the right-most integral above
by using the expansion in \eqref{eq:psi} as follows:
$$
\iint\limits_{|\zeta|>1}|\psi(\zeta)|^2\,d\xi d\eta
= \pi\sum_{n=1}^{\infty}n|a_n|^2.
$$
Plugging this with \eqref{eq:Chichra}, \eqref{eq:area} and \eqref{eq2d},
we obtain
$$
(k^{-2}-1)\pi\sum_{n=1}^{\infty}n|a_n|^2\le
\pi\left[ \frac{1}{(1-p^2)^2}-\sum_{n=1}^{\infty}n|a_n|^2\right],
$$
which yields the desired inequality.

Finally, we analyze the equality case for \eqref{eq2a0}.
Suppose that equality holds in \eqref{eq2a0}.
Then, equalities must hold both in \eqref{eq:area} and in \eqref{eq2d}.
The equality in \eqref{eq2d} implies that $\dz\psi=0$ on $\ID.$
In other words, $h=\bar\psi$ is analytic on $\ID.$
Therefore, $\phi(\zeta)=\zeta/(1-p\zeta)+\overline{h(\zeta)}.$
The equality in \eqref{eq:area} means that $|\dzb\phi/\dz\phi|$ is
the constant $k$ a.e.~on $\ID.$
Since $\dzb\phi(\zeta)/\dz\phi(\zeta)=\overline{h'(\zeta)}(1-p\zeta)^2,$
it implies that the analytic function $h'(\zeta)(1-p\zeta)^2$ has constant
modulus $k$ and therefore a constant $\alpha$ with $|\alpha|=k.$
Hence, $h'(\zeta)=\alpha(1-p\zeta)^{-2}$ for $|\zeta|<1.$
Integrating it, we obtain $h(\zeta)=\alpha\zeta/(1-p\zeta)+h(0).$
Thus, we finally have the form
$$
\phi(\zeta)=\frac{\zeta}{1-p\zeta}+\overline{h(0)}
+\frac{\bar\alpha\bar\zeta}{1-p\bar\zeta},
\quad \zeta\in\ID.
$$
Therefore,
$$
f(z)=\frac{1}{z-p}+\overline{h(0)}+\frac{\bar\alpha}{\bar z-p},
\quad z\in\ID^*,
$$
whose boundary values on $\partial\ID$
are the same as those of the meromorphic function
$$
g(z)=\frac{1}{z-p}+\overline{h(0)}+\frac{\bar\alpha z}{1-pz}.
$$
Since $f(z)-g(z)$ is bounded analytic on $\ID,$
$f(z)$ is identically equal to $g(z)$ on $\ID$ by the maximum principle.
In particular, $h(0)=\overline{a_0}$ and $\alpha=\overline{a_1}.$
Thus we have seen that the function $f$ must have the form \eqref{eq2a1}
if equality holds in \eqref{eq2d}.
We need to show that the function $f$ of the form \eqref{eq2a1}
is indeed a member of $\Sigma_k(p).$
We first show that $f$ is univalent in $\ID.$
We compute
$$
f(z_1)-f(z_2)=\frac{z_2-z_1}{(z_1-p)(z_2-p)}
\left[1-a_1\left(\frac{z_1-p}{1-pz_1}\cdot\frac{z_2-p}{1-pz_2}\right)\right]
$$
for $z_1,z_2\in\ID.$
Since
$$
\left|a_1\left(\frac{z_1-p}{1-pz_1}\cdot\frac{z_2-p}{1-pz_2}\right)\right|
<|a_1|=k\le 1
$$
for $z_1, z_2\in\ID,$ we see that $f(z_1)\ne f(z_2)$ if $z_1\ne z_2.$
Hence, $f\in\Sigma(p).$
On the other hand, the function in \eqref{eq:ext} agrees with
that in \eqref{eq2a1} on the boundary $|z|=1,$ and
is a composition of the M\"obius transformation
$1/(z-p)$ with the $k$-quasiconformal affine mapping $w+a_0+a_1\bar w.$
Thus we conclude that $f$ belongs to $\Sigma_k(p).$
\epf

\section{Proof of Theorems \ref{Theorem-2} and \ref{Theorem-3}
and a concluding remark}

We start with the proof of Theorem \ref{Theorem-2}.

\bpf[Proof of Theorem \ref{Theorem-2}]
We first investigate the behavior of the function $\omega.$
Since $|\omega'|\le k(1+p)^{-2},$ we have the inequality
$$
|\omega(z_1)-\omega(z_2)|\le \frac{k}{(1+p)^2}|z_1-z_2|,
\quad z_1,z_2\in\ID.
$$
We thus conclude that $\omega$ is (uniformly) Lipschitz continuous on $\ID.$
In particular, $\omega$ extends to $\overline\ID$ continuously with
the same Lipschitz constant.
Also, the function $f(z)=1/(z-p)+\omega(z)$ extends to $|z|\le1$ continuously.
We now see that
$$
\left|\frac1{z_1-p}-\frac1{z_2-p}\right|
= \frac{|z_1-z_2|}{|z_1-p||z_2-p|}\ge \frac{|z_1-z_2|}{(1+p)^2},
\quad z_1,z_2\in\ID.
$$
Hence, we have
$$
|f(z_1)-f(z_2)|\ge \frac{1-k}{(1+p)^2}|z_1-z_2|
$$
for $z_1,z_2\in\overline\ID,$ which proves that $f(z)$ is univalent on
$|z|\le 1.$

For a while, we denote by $g(z)$ the function appearing in the right-hand
side of \eqref{eq0}.
Obviously, $g(z)$ can be defined in $|z|\ge1$ and agrees
with the above $f(z)$ on $|z|=1.$
We show now that $g$ is $k$-quasiconformal on $|z|>1.$
A computation yields
$$
\dz g(z)=\frac{-1}{(z-p)^2}
\aand
\dzb g(z)=\frac{-\omega'(1/\bar z)}{\bar z^2}.
$$
Hence, the complex dilatation $\mu=\dzb g/\dz g$ of $g$ satisfies
$$
|\mu(z)|=\frac{|z-p|^2}{|z|^2}|\omega'(1/\bar z)|
\leq (1+p)^2|\omega'(1/\bar z)|
\le k <1
$$
for $|z|>1.$
Since the unit circle is removable for quasiconformality
(see \cite[p.205]{LV:qc}),
we conclude that $g$ gives a $k$-quasiconformal extension of $f.$
\epf

A straightforward application of Theorem \ref{Theorem-2} yields the following sufficient condition for a function $f$ of the form (\ref{eq0a}) to
belong to $\Sigma_k(p)$.

\bcor\label{cor-1}
Let $0\le p<1$ and $0\le k<1.$
Suppose that a meromorphic function $f(z)$ on $|z|<1$ has the form
\eqref{eq0a}.
If
$$
\sum_{n=1}^{\infty}n|a_n|\leq \frac{k}{(1+p)^2},
$$
then $f\in \Sigma_k(p)$.
\ecor

\bpf
This immediately follows from Theorem \ref{Theorem-2} because
$$
|\omega'(z)|
\leq \sum_{n=1}^{\infty}n|a_n||z|^{n-1}
\leq \sum_{n=1}^{\infty}n|a_n|\le \frac{k}{(1+p)^2},\quad |z|<1.
$$
\epf

Next we prove Theorem \ref{Theorem-3}.

\bpf[Proof of Theorem \ref{Theorem-3}]
Let $f\in \Sigma_{k_1}(p)$ and $g\in \Sigma_{k_2}(p)$ be expressed
as in \eqref{eq2a1}.
Then Theorem~\ref{Theorem-1} gives us
$$
\sum_{n=1}^{\infty}n|a_n|^2\leq \frac{{k_1}^2}{(1-p^2)^2}
\aand
\sum_{n=1}^{\infty}n|b_n|^2\leq \frac{{k_2}^2}{(1-p^2)^2}.
$$
Now an application of Cauchy-Schwarz inequality together with the
aforementioned inequalities yields
\beqq
\sum_{n=1}^{\infty}n|a_n b_n|
&=&\sum_{n=1}^{\infty}(\sqrt{n}|a_n|)(\sqrt{n}|b_n|)\\
&\leq& \left(\sum_{n=1}^{\infty} n|a_n|^2\right)^{1/2}
\left(\sum_{n=1}^{\infty} n|b_n|^2\right)^{1/2}\\
&\leq& \frac{k_1k_2}{(1-p^2)^2}
= \frac{\alpha}{(1+p)^2},
\eeqq
where $\alpha=k_1k_2(1-p)^{-2}$.
Since $\alpha<1$ by assumption,
the desired result follows from Corollary \ref{cor-1}.
\epf

We conclude the present note with an outline of another proof of
\eqref{eq2a0} based on Lehto's principle
(cf.~\cite[II.3.3]{Lehto:univ}) and Theorem B.
Before it, we recall the definition of the complex Banach
(indeed, Hilbert) space $\ell^2.$
This is the set of sequences $x=\{x_n\}_{n=1}^\infty$
of complex numbers with the norm
$$
\|x\|_{\ell^2}=\left(\sum_{n=1}^\infty |x_n|^2\right)^{1/2}
<\infty.
$$

It is enough to show \eqref{eq2a0} for functions in
$\Sigma_k^0(p)$ only.
Suppose that $f\in\Sigma_k^0(p)$ is already extended to a
$k$-quasiconformal mapping of $\sphere$
and let $\mu$ be its complex dilatation.
We remark that $|\mu|\le k$ a.e.\,in $\ID^*$ and $\mu=0$
in $\ID.$
By the measurable Riemann mapping theorem, for each $t\in\ID,$
there exists a unique quasiconformal mapping $f_t$ of $\sphere$
for which the complex dilatation is $t\mu/k$ and $f_t|_{\ID}
\in\Sigma^0(p).$
Note here that $f_k=f.$
Then $f_t$ has an expansion of the form
$$
f_t(z)=\frac1{z-p}+\sum_{n=1}^\infty a_n(t)z^n
$$
in $|z|<1.$
By the holomorphic dependence of the solution to the Beltrami
equation, $a_n(t)$ is analytic in $|t|<1$ for every $n\ge1.$
We now consider the sequence
$\sigma(t)=\{\sqrt n a_n(t)\}_{n=1}^\infty.$
Theorem B tells us that $\|\sigma(t)\|_{\ell^2}\le 1/(1-p^2).$
Hence, we conclude that $\sigma:\ID\to\ell^2$
is a bounded analytic function taking values in the complex
Banach space $\ell^2.$
Since $\sigma(0)=0,$ the (generalized) Schwarz lemma
yields the inequality $\|\sigma(t)\|_{\ell^2}\le |t|/(1-p^2).$
In particular, letting $t=k$ gives \eqref{eq2a0}.

We must say that this method is conceptually simpler than
that of our proof in Section 2.
However, this does not provide information about the
equality case in an obvious manner.

\bigskip

{\bf Acknowledgement:} The authors thank Karl-Joachim Wirths for his
suggestions and careful reading of the manuscript.


\def\cprime{$'$} \def\cprime{$'$} \def\cprime{$'$}
\providecommand{\bysame}{\leavevmode\hbox to3em{\hrulefill}\thinspace}
\providecommand{\MR}{\relax\ifhmode\unskip\space\fi MR }
\providecommand{\MRhref}[2]{%
  \href{http://www.ams.org/mathscinet-getitem?mr=#1}{#2}
}
\providecommand{\href}[2]{#2}

\end{document}